\documentclass{amsart}
\usepackage{amsmath,amsthm}
\usepackage{amsfonts,amssymb}

\newtheorem{thm}{Theorem}[section]

\newtheorem{lem}[thm]{Lemma}
\newtheorem{con}[thm]{Conjecture}
\newtheorem{prop}[thm]{Proposition}

\newtheorem{defn}[thm]{Definition}

\theoremstyle{remark}

 \def\CS{{\mathcal S}}

 \def\NN{{\mathbb N}}

 \def\RR{{\mathbb R}}
 \def\ZZ{{\mathbb Z}}

        \def\sign{\operatorname{sign}}

\begin{document}
\input epsf    
To appear in J. Experimental Math. 
\medskip

\title[]
{On polynomials of least deviation from zero in several variables}
\author{Yuan Xu}
\address{Department of Mathematics\\ University of Oregon\\
    Eugene, Oregon 97403-1222.}\email{yuan@math.uoregon.edu}

\date{Feburary 25, 2003}
\thanks{Work partially supported by the National Science Foundation under 
Grant DMS-0201669}

\begin{abstract}
A polynomial of the form $x^\alpha - p(x)$, where the degree of $p$ is less
than the total degree of $x^\alpha$, is said to be least deviation from
zero if it has the smallest uniform norm among all such polynomials. We study 
polynomials of least deviation from zero over the unit ball, the unit sphere 
and the standard simplex. For $d=3$, extremal polynomial for $(x_1x_2x_3)^k$ 
on the ball and the sphere is found for $k=2$ and $4$. For $d \ge 3$, a family
of polynomials of the form $(x_1\cdots x_d)^2 - p(x)$ is explicit given and 
proved to be the least deviation from zero for $d =3,4,5$, and it is 
conjectured to be the least deviation for all $d$. 
\end{abstract}

\maketitle 

\section{Introduction} 
\setcounter{equation}{0} 

Let $\Pi_n^d$ denote the space of polynomials of degree at most $n$ in $d$ 
variables and we write $\Pi_n = \Pi_n^1$. For $d =1$, it is well-known that 
the $2^{1-n}$ multiple of the Chebyshev polynomial of the first kind 
$$
   T_n(x) = \cos n(\arccos x) = 2^{n-1} x^n + q(x), \qquad q \in \Pi_{n-1},
$$
is the monic polynomial of least deviation from zero in $\Pi_n$ in the 
space $C[-1,1]$; that is, 
$$
  \inf_{p\in \Pi_{n-1}} \|x^n - p(x)\|_{C[-1,1]} = 2^{1-n} \|T_n\|_{C[-1,1]}
    = 2^{1-n}.
$$
Equivalently, we say that $x^n - 2^{1-n}T_n$ is the best approximation to $x^n$
in $C[-1,1]$.

Let $\Omega$ be a region in $\RR^d$. For $f \in C(\Omega)$, the best 
approximation of $f$ from $\Pi_n^d$ in the uniform norm is the quantity 
\begin{align} \label{eq:1.1} 
 E_n(f;\Omega) = \inf_{p \in \Pi_{n-1}^d} \|f - p\|_{C(\Omega)}, 
\end{align} 
where $\|f\|_{C(\Omega)} = \max_{x \in \Omega} |f(x)|$. We call $p^*$ an 
extremal polynomial for $f$ if $E_n(f;\Omega)=\|f-p^*\|_{C(\Omega)}$.
For $x = (x_1,\ldots, x_d) \in \RR^d$ and $\alpha =(\alpha_1,\ldots,\alpha_d)
\in \NN_0^d$, we define the monomial $x^\alpha = x_1^{\alpha_1}\cdots 
x_d^{\alpha_d}$. The degree of the monomial $x^\alpha$ is $|\alpha| = 
\alpha_1 +\ldots+ \alpha_d$. 
If $p^*(x)$ is an extremal polynomial for the monomial $x^\alpha$, we call 
$x^\alpha - p^*(x)$ the polynomial of least deviation from zero. For $\Omega$ 
being a region in $\RR^d$, polynomials of least deviation are known only 
in the case that $\Omega$ is a cube. We are interested in the case of the 
unit ball $B^d =\{x: \|x\| \le 1\}$, where $\|x\|$ is the usual Euclidean norm
of $x$, the unit sphere $S^{d-1} =\{x : \|x\| = 1\}$ and the standard simplex 
$T^d =\{x: x_1 \ge 0, \ldots, x_d\ge 0, 1-x_1-\ldots -x_d \ge0\}$. 

For $d=2$, the least deviation of $x^n y^m$ from $\Pi_{n+m-1}^2$ in the space 
$C(\Omega)$ has been studied for $B^2$ and $T^2$ (see, for example, 
\cite{G}, \cite{R}, \allowbreak \cite{NX}, \cite{BHN}). For $d>2$, the only case known is 
$x_1\cdots x_d$ on $B^d$ and $S^{d-1}$, which is a polynomial of least 
deviation by itself. This is shown recently in \cite{AY},
\begin{align} \label{eq:1.2} 
\inf_{p\in \Pi_{d-1}^d }\|x_1 \cdots x_d - p(x)\|_{C(B^d)} & = 
\inf_{p\in \Pi_{d-1}^d }  \|x_1 \cdots x_d - p(x)\|_{C(S^{d-1})} \\
 & = \|x_1 \cdots x_d\|_{C(S^{d-1})} =d^{-d/2}. \notag
\end{align}
In other word, the best approximation of $x_1 \cdots x_d$ from $\Pi_{d-1}^d$ 
is the zero polynomial. Finding polynomials of least deviation on these regions
appear to be a difficult problem. Only a handful of explicit non-trivial 
examples of extremal polynomials for $d \ge 3$ are known in the literature.

In the present paper, we study the least deviation from zero for monomials of
lower degrees. We found extremal polynomials for $(x_1 x_2 x_3)^2$ 
and $(x_1 x_2 x_3)^4$ on $B^3$ and $S^2$ and a family of extremal polynomials 
for $x_1^2\cdots x_d^2$ on $B^d$ and $S^{d-1}$, which are derived from the 
extremal polynomials for $x_1x_2x_3$ and $(x_1x_2x_3)^2$ on $T^3$ and 
$x_1\cdots x_d$ on $T^d$, respectively. We give an explicit construction of 
this family of polynomials and conjecture that 
they are the least deviation polynomials. The conjecture is proved for $d =3, 
4, 5$. The result provides, we believe, the first non-trivial example of 
polynomials of least deviation on these domains. For example, we have
$$
\inf_{p \in \Pi_5^3} \|x_1^2 x_2^2 x_3^2 - p(x) \|_{C(B^3)}= 
  \inf_{p \in \Pi_5^3} \|x_1^2 x_2^2 x_3^2 - p(x) \|_{C(S^3)}
   = 72^{-2}
$$
and the minimum is attained by the extremal polynomial $R_3(x,y)$ defined by  
$$
 R_3(x_1,x_2,x_3) = 72 x_1^2 x_2^2 x_3^2 - 4(x_1^2+x_2^2+x_3^2) +
   4(x_1^4+x_2^4+x_3^4)^2 + 1. 
$$
The least deviation, $72^{-2}$, is surprisingly small in view of the value 
$3^{3/2}$ for $x_1 x_2 x_3$. 

Our proof is based on a general result for the Chebyshev approximation in 
\cite{RS}, in which the best approximation element is characterized in terms 
of extremal signature. The most difficult part, however, is to identify a 
correct extremal polynomial. There is no general method for this purpose. 
We relied heavily on the computer algebra system {\it Mathematica} to test 
and verify conjectures. In retrospect, the explicit construction is natural 
and rather suggestive. For example, $R_3(x)$ agrees with the Chebyshev 
polynomial $T_2(x)$ on the three edges of the face of $T^3$ defined by 
$x_1+x_2+x_3 =1$. The result allows first glimpse of what an extremal 
polynomial in more than two variables may look like. 

The paper is organized as follows. In the following section we recall the
theoretic background needed to prove our result. The results for $d=3$ are
discussed in Section 3 and those for $d >3$ are in Section 4. 

\section{Extremal signature and best approximation}
\setcounter{equation}{0} 

We recall the characterization of the extremal polynomials in terms of the
extremal signature. The study in \cite{RS} is given in the general setting
of approximation from a finite dimensional subspace of $C(\Omega)$ on a 
compact Hausdorff space $\Omega$. We shall restrict the statement to our
setting. 

Let $\Omega$ be an infinite compact set in $\RR^d$. A signature $\sigma$ on 
the set $\Omega$ is a function with finite support, whose nonzero values are 
either $+1$ or $-1$. A signature $\sigma$ is called extremal with respect to 
$\Pi_n^d$ if there exists a subset $\CS$ in the support of $\sigma$ and 
positive numbers $\lambda_v$, $v\in \CS$, such that 
$$
   \sum_{v \in \CS} \lambda_v \sigma(v) p(v) = 0, \qquad 
     \hbox{ for all $p$ in $\Pi_n^d$}. 
$$
Let $r >0$ be a fixed number. For each $p \in \Pi_{n-1}^d$, we denote by
$\CS_r(p;f)$ the set
$$
  \CS_r (p;f) = \{x \in \Omega: |f(x) - p(x)| = r\}. 
$$
If $r = \|f - p\|_{C(\Omega)}$, $\CS_r(p;f)$ is the set of extremal points of 
$f -p$ and we denote it by $\CS(p;f)$. 

The characterization of the best approximation of $f$ from $\Pi_n^d$ is
given by the following theorem in \cite{RS}:

\begin{thm} \label{thm:2.1}
A polynomial $p^*$ in $\Pi_n^d$ satisfies $\|f-p^*\|_{C(\Omega)} = 
E_n(f; \Omega)$ if and only if there exists an extremal signature $\sigma$ 
with support in $\CS(p^*;f)$ such that $\sigma(v) = \sign (f-p^*)(v)$ for 
all $v \in \CS(p^*;f)$. 
\end{thm}

The sufficient part of the theorem provides a method to verify if a polynomial
$p^*$ is extremal. One needs, however, to know the extremal polynomial in 
advance, as the extremal signature is supported on the set $\CS(p^*;f)$ which 
depends on $p^*$. The sufficient part of the theorem can be extended to the 
signature support on $\CS_r(p;f)$, in which $r$ is not necessarily 
$\|f-p\|_{C(\Omega)}$. We will use this slightly extended version, which we
state in the following. A simple proof is included for completeness; see 
\cite{RS} for more details.

\begin{thm} \label{thm:2.2}
Suppose there exists a polynomial $p^* \in \Pi_n^d$ and an extremal signature
$\sigma$ supported on $\CS_r(p^*;f)$. Then $E_n(f; \Omega)\ge r$.
\end{thm}

\begin{proof}
We can normalize the measure $\lambda_\mu$ for the extremal signature so that 
it is a probability measure; that is, $\sum_{v\in S_r(p^*;f)} \lambda_v =1$. 
Let $S(r) = S_r(p^*,f)$ in this proof. Since $\sum \lambda_v p(v) =0$ for any 
polynomial $p \in \Pi_n^d$, we have 
\begin{align*}
\|f(x) - p(x)\|_{C(\Omega)} & \ge \sum_{v \in \CS(r)} \lambda_v |f(v) - p(v)|\\
  &\ge \Big| \sum_{v \in \CS(r)} \lambda_v \sigma_v f(v) - 
  \sum_{v \in \CS(r)} \lambda_v \sigma_v p(v) \Big| \\
 & =  \Big| \sum_{v \in \CS(r)} \lambda_v \sigma_v f(v)\Big| = 
  \Big| \sum_{v \in \CS(r)} \lambda_v \sigma_v (f(v) - p^*(v)) \Big|\\
 & = \sum_{v \in \CS(r)} \lambda_v |f(v) - p^*(v)| = r \sum_{v \in \CS(r)} 
   \lambda_v = r. 
\end{align*}
where we have used the fact that $f(v) - p^*(v) = \sigma_v r$ for $v\in S(r)$.
\end{proof}

The extension allows us to apply the result to the situation where a good 
candidate for $p^*$ is identified but the norm of $f-p^*$ is hard to determine.
This is precisely our case in Section 4. 

Our construction is motivated by the recent study in \cite{AY}, in which it is
shown that if $f$ is invariant under a finite group $G$ (that is, $f(x g) = 
f(x)$ for all $g \in G$), then the best approximation $E_n(f;S^{d-1})$ is 
attained at $G$ invariant polynomials. (This result appeared early in \cite{GP}
as pointed out by a referee.) More precise, we state the result in
\cite{AY} as follows:

\begin{prop} \label{prop:2.3}
Let $G$ be a subgroup of the rotation group $O(d)$ and let $G\Pi_n^d$ denote 
the polynomials in $\Pi_n^d$ that are invariant under $G$. If $f$ is invariant
under $G$, then
$$
\inf_{p \in \Pi_{n-1}^d} \|f(x) - p(x)\|_{C(S^{d-1})} =
  \inf_{p \in G\Pi_{n-1}^d} \|f(x) - p(x)\|_{C(S^{d-1})}. 
$$
\end{prop}

Using this fact, the best approximation of several invariant functions are
given in \cite{AY}, including the case $x_1 \cdots x_d$ in \eqref{eq:1.2} 
(invariant under the symmetric group). The proof in \cite{AY} can be applied to
any region $\Omega$ and $f$ that are invariant under a finite group $G$. In 
particular, if $f$ is invariant under a subgroup $G$ of the symmetric group 
$S_d(T^d)$ of the simplex $T^d$, then an extremal polynomial of $f$ can be 
taken as a $G$-invariant polynomial. 

If $f$ is even in each of its variables, then $f$ is invariant under the 
sign changes of each variable (invariant under the group $\ZZ_2^d$); the 
extremal polynomial can be taken as a polynomial even in each of its 
variables. Furthermore, instead of $B^d$ or $S^{d-1}$, we can work with $T^d$
and $T^{d-1}$ in this case. In fact, the following general proposition holds: 

\begin{prop} \label{prop:2.4}
Let $\alpha \in \NN_0^d$ and write $2 \alpha = (2\alpha_1, \ldots,2\alpha_d)$
and $|\alpha|=n$. If $p^*(x)$ is an extremal polynomial for $E_n(x^\alpha;T^d)$
then $p^*(x_1^2,\ldots,x_d^2)$ is an extremal polynomial for 
$E_n(x^{2\alpha}; B^d)$; conversely, if $q^*$ is an extremal polynomial for 
$E_n(x^{2\alpha};B^d)$ in the form $q^*(x) = p^*(x_1^2,\ldots,x_d^2)$, then
$p^*(x)$ is an extremal polynomial for $E_n(x^\alpha; T^d)$. 
Furthermore, let $f_{\alpha}(x_1,\ldots,x_{d-1}) = x_1^{\alpha_1}\cdots 
x_{d-1}^{\alpha_{d-1}} (1-x_1-\ldots-x_{d-1})^{\alpha_d}$; then the above 
conclusion holds for $E_n(f_{\alpha}; T^{d-1})$ and $E_n(x^{2\alpha};S^{d-1})$.
\end{prop}

We note that $f_\alpha(x_1^2, \ldots,x_{d-1}^2) = x^{2 \alpha}$ on $S^{d-1}$.  
The proposition follows easily from the fact that $x \mapsto (x_1^2,\ldots,
x_d^2)$ is one-to-one from $T^d$ to $B_+^d = \{x \in B^d: x_i \ge0\}$, and 
the map also induces a one-to-one mapping from $\Pi_n^d$ to $G\Pi_{2n}^d$ with
$G= \ZZ_2^d$, that is, the subspace of polynomials that are even in each of its
variables. For $d =2$ the proposition has been used in \cite{BHN}. The 
correspondence between polynomials on these domains also works for other 
problems involving polynomials, such as orthogonal polynomials and cubature
formulae, see for example \cite{X98b}.

\section{Least deviation from zero for $d =3$}
\setcounter{equation}{0} 

We consider best approximation to the monomials $x_1x_2x_3$ and 
$(x_1x_2x_3)^2$ in this section. The main task is to identify an extremal
polynomial. The results in the previous section provide some guidance, 
but there is no general method for this purpose. Our first example, $R_3(x)$ 
given below, was found after many attempts. See the comments after the proof.

\begin{thm} \label{thm:3.1}
Define the polynomial $R_3(x)$ by
$$
 R_3(x) =  72 x_1 x_2 x_3 - 4 (x_1+x_2+x_3) + 4 (x_1+x_2+x_3)^2  
 - 8 (x_1 x_2+ x_2 x_3+ x_1x_3) + 1.  
$$
Then $72^{-1} R_3(x)$ is a polynomial of least deviation from zero and
$$
E_2(x_1x_2x_3; T^3) = E_2(x_1x_2(1-x_1-x_2); T^2) =72^{-1}\|R_3\|_{C(T^3)}
 = 72^{-1}
$$
Furthermore, $72^{-2}R_3(x_1^2,x_2^2,x_3^2)$ is a polynomial of least 
deviation from zero and  
$$
E_5(x_1^2x_2^2x_3^2; B^3) = E_5(x_1^2x_2^2x_3^2; S^2) = 72^{-2}
 \left\|R_3(x_1^2,x_2^2,x_3^2) \right\|_{C(B^3)} = 72^{-2}.
$$
\end{thm}

\begin{proof}
By Proposition \ref{prop:2.4} we only need to work with the simplex. It is 
easy to verify that $R_3(0,0,0)=1$ and $R_3(1/2,1/2,0) = -1$. Solving the 
equations $\partial_i R_3(x) =0$, $i =1,2,3$, shows that $R_3$ has 4 critical
points inside $T_3$ but none of them is maximum or minimum, since the values 
of $|R_3(x)|$ at these points are less than 1. Thus, $|R_3(x)|$ attains its
maximum on the boundary of $T^3$. It is easy to verify that the polynomial 
$R_3(x)$ satisfies 
$$
  R_3(x,y,0) = R_3(x,0,y) = R_3(0,x,y) = (1 - 2x)^2 + (1-2y)^2 - 1
$$
which is bounded by 1 in absolute value. Hence, we only need to show
that $|R_3(x)|$ is bounded by one on the face of $T^3$ defined by
$x_1+x_2+x_3 = 1$; that is, we need to show that $U_3(x_1,x_2) =
R_3(x_1,x_2,1-x_1-x_2)$ is bounded by 1 in absolute value on $T^2$.
Taking derivatives of $U_3(x_1,x_2)$ and solving for the critical 
points shows that it has 4 critical points inside $T^2$, of which 
only the point $(1/3,1/3)$ is a maximal, $U(1/3,1/3)=1$. Furthermore,
it is easy to verify that 
$$
 U_3(x,0) = U_3(0,x) = U_3(x,1-x) = T_2(x); 
$$
that is, it agrees with Chebyshev polynomial of degree $2$ on the 
boundary of the triangle. Hence, $|U_3(x)| \le 1$. Furthermore, the
above analysis also shows that 
\begin{align*}
\CS_+& := \{x : |R_3(x) =1\} = \{(0,0,0),(1,0,0),(0,1,0),(0,0,1),
(1/3,1/3,1/3)\}\\
\CS_-& := \{x : |R_3(x) =-1\} = \{(1/2,1/2,0),(1/2,0,1/2),(0,1/2,1/2)\}.
\end{align*}  
Let $\sigma(v) = 1$ on $\CS_+ \setminus\{(0,0,0)\}$ and $\sigma(v) = -1$ on 
$\CS_-$. We show that $\sigma$ is an extremal signature. Define $L_1f$ and 
$L_2 f$ by 
$$
 L_1 f = \frac{3}{4} f\left(\frac{1}{3},\frac{1}{3},\frac{1}{3}\right) +
   \frac{1}{12} (f(1,0,0) + f(0,1,0) + f(0,0,1))
$$
and 
$$ 
L_2 f = \frac{1}{3}\left( f\left(\frac{1}{2},\frac{1}{2},0\right)+ 
    f\left(\frac{1}{2},0,\frac{1}{2}\right) + 
    f \left(0,\frac{1}{2},\frac{1}{2}\right)\right). 
$$
Then $L f : = L_1 f - L_2 f$ satisfies $L f =0$, $f \in \Pi_2^3$. Thus, 
$\sigma$ is an extremal signature for $x_1x_2x_3$ in $C(T^3)$. Furthermore, 
the support sets $\CS_+$ and $\CS_-$ of $\sigma$ are on the face of $T^3$, 
which is identified with $T^2$. This shows that $\sigma$ is also an extremal 
signature for $x_1x_2(1-x_1-x_2)$ in $C(T^2)$. 
\end{proof}

Let us mention a connection between cubature formulae and the extremal 
signature for $R_3(x)$. A cubature formula is a linear combination of function 
evaluations that gives an approximation to an integral (\cite{S}). Let $d\mu$ 
be a positive measure on $\Omega \subset \RR^d$. If 
$$
\int_{\Omega} f(x) d\mu  = \sum_{k=1}^N \lambda_k f(x_k), \qquad f \in \Pi_n^d,
$$
and there is at least one $f \in \Pi_{n+1}^d$ such that the equality fails,
then the cubature formula is said to be of degree $n$. It is called positive,
if all $\lambda_k$ are positive numbers. For the extremal signature for 
$R_3(x)$, it is easy to verify that both $L_1f$ and $L_2f$ are cubature 
formulae of degree $2$ for $d x$ on the set $\Sigma^2 = \{x \in T^3: 
x_1+x_2+x_3 =1\}$; that is, 
$$
\int_{\Sigma^2} f(x)dx = L_1 f = L_2f, \qquad \hbox{for all $f$ in $\Pi_2^2$}.
$$
Since we identify $\Sigma^2$ with $T^2$, one can write $L_1$ and $L_2$ as 
linear combinations of function evaluations for functions of two variables. 
Thus, the extremal signature is given by the difference of two positive 
cubature formulae.

During our search for $R_3(x)$, we found $U_3(x)$ first. In retrospect, the
formula of $U_3(x)$, which can be written as 
$$
U_3(x_1,x_2) =  72 x_1 x_2 (1-x_1-x_2) - 3 + 4 (x_1^2+ x_2^2+(1-x_1-x_2)^2),  
$$
is quite natural since it agrees with the Chebyshev polynomials of degree $2$ 
on the boundary of $T^2$. This also suggests the possibility that other 
monomials may also have extremal polynomials that agree with Chebyshev 
polynomials on the boundary of $T^3$. For example, for $(x_1 x_2 x_3)^n$, one
may look for a polynomial that agrees with Chebyshev polynomials of $n$-th 
degree on the boundary of the simplex. One example of such a polynomial is 
$T_n(R_3(x))$, where $T_n(t)$ denotes the Chebyshev polynomial of degree $n$.
Although this function is not a polynomial of least deviation, it helps
us to find a solution for the monomial $(x_1 x_2 x_3)^2$. 

\begin{thm} \label{thm:3.2}
Define polynomial $R_5(x)$ by 
\begin{align*}
R_5(x) = & 27^2 b\, (x_1 x_2 x_3)^2 -1+2(x_1+ x_2 + x_3)-2(x_1+ x_2 + x_3)^2\\
 & + 2\left[1 - 4(x_1 + x_2 + x_3)+4(x_1^2 + x_2^2 + x_3^2)\right]^2 \\ 
 &- 27 x_1 x_2 x_3\left[ 
 (32/9 - 2 a + b)(x_1 + x_2 + x_3)^2 + 6a(x_1 x_2 + x_1 x_3 + x_2 x_3)\right].
\end{align*}
Then $27^{-2}b^{-1} R_5(x)$ is a polynomial of least deviation from zero,
$$
E_5(x_1^2x_2^2x_3^2; T^3) = E_5(x_1^2x_2^2(1-x_1-x_2)^2; T^2) =
 27^{-2}b^{-1} \|R_5\|_{C(T^3)} = 27^{-2}b^{-1}, 
$$
where the constant $a, b$ and the reciprocal of the least deviation is given by
$$
a = 28.5926243,\qquad b=21.8935834, \qquad 27^2 b =15960.4223.
$$
Furthermore, $72^{-4}b^{-2}R_3(x_1^2,x_2^2,x_3^2)$ is a polynomial of least 
deviation from zero,
$$
E_5(x_1^2x_2^2x_3^2; B^3) = E_5(x_1^2x_2^2x_3^2; S^2) = 72^{-4}b^{-2}
 \left\|R_5(x_1^2,x_2^2,x_3^2) \right\|_{C(B^3)} = 27^{-4}b^{-2}.
$$
\end{thm}

Just like the case of $R_3(x)$, the proof amounts to showing that 
$|R_5(x)|\le 1$ on $T^3$ and there exists an extremal signature. It is not 
difficult once the formula of $R_5$ is identified. We first give an account 
on how $R_5$ is discovered.  

Following the construction of $R_3(x)$, we look for a polynomial in the form of
\begin{align*}
 U_5(x,y) = & 27xy(1-x-y )\left[27 b\, x y (1 - x - y)+ 
   3a(x^2+y^2 +(1-x-y)^2)- c \right]\\
   &  + 2 \left[ -3 + 4 (x^2 + y^2 + (1-x-y)^2)\right]^2 -1     
\end{align*}
that will be a polynomial of least deviation on $T^2$ with leading monomial 
$x^2 y^2 (1-x-y)^2$. Note that the Chebyshev polynomial of degree $2$ is 
$T_2(t) = 2t^2 -1$ and $T_2(2t-1) = -3+4(t^2 + (1-t)^2))$. The form of $U_5$ 
is chosen so that on the boundary of $T^2$ it satisfies 
$$
 U_5(x,0) = U_5(0,x) = U_5(x,1-x) = T_2(T_2(2x-1)) = T_4(2x-1). 
$$
We then choose $c = 2/9 + a + b$ so that $U_5(1/3,1/3) =1$. It follows that 
$1-U_5(x,x)$ can be factored as
$$
  1-U_5(x,x) = x(1-2x)(1- 3 x)^2(64 - 54 a x + 27 b x + 162 b x^2).
$$
We need to choose $a$ and $b$ so that the last factor is positive for 
$0 \le x \le 1/2$. One choice is to make this factor $2b(9 x + d)^2$. This 
leads to $a= 16(3-4d)/(3d^2)$ and $b = 32/d^2$. At this point, it becomes
apparent that there need to be more points on which $U_5(x,y) = -1$ inside
$T^2$. We therefore solve the equations $U_5(x,x) = -1$ and $U_5'(x,x) =0$.
This leads to $d = -1.208972894$, which gives the values for $a$ and $b$ in 
the theorem. It turns out that this choice does work out and $|U_5(x,y)|\le 1$ 
on $T^2$. The final step is to identify the formula of $R_5(x_1,x_2,x_3)$ from
that of $U_5(x_1,x_2)$ with the requirement that $R_5(x_1,x_2,1-x_1-x_2)=
U_5(x_1,x_2)$ and $|R_5(x,y)| \le 1$ on $T^3$. This step is not trivial since 
an additional multiple of $(x_1+x_2+x_3)^k$ to any term in $R_5(x)$ does not
change the value of the polynomial on the face of $T^3$ defined by $x_3 =
1-x_1-x_2$. The polynomial $U_5(x_1,x_2)$ is an extremal polynomial on the
triangle $T^2$ that agrees with the Chebyshev polynomials of degree 4 on 
the three boundary lines of $T^2$, its graph is depicted below.

\smallskip
\centerline{
\epsfxsize= 4in
\epsffile{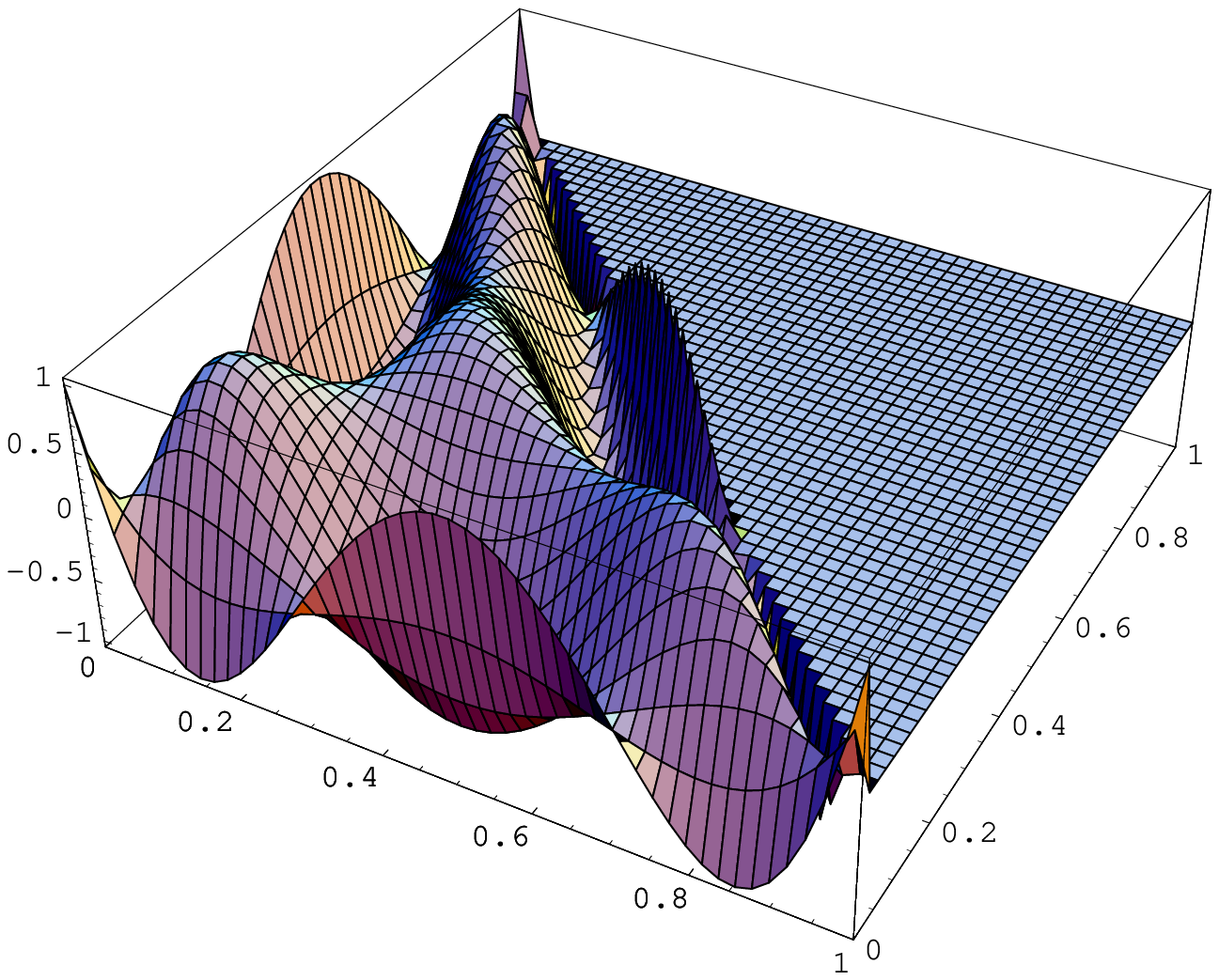} }
\centerline{ The polynomial $U_3$}
\smallskip

Let us point out that there does not seem to be a closed form for the values 
of $a$ and $b$. In fact, the value of $d$ in the above paragraph is one of 
the real roots of the following polynomial: 
\begin{align*}
& -612220032 - 1365527808\, t - 835528041 \, t^2 - 101556504 \, t^3 \\
& +    23270976 \, t^4  
 + 26037504 \, t^5 + 7670016 \, t^6 + 929280 \, t^7 +  41984 \, t^8.
\end{align*}
This polynomial has 4 real roots and 4 complex roots, and it cannot be 
factored over the integers.

We now give a formal proof of Theorem \ref{thm:3.2}. 

\medskip\noindent
{\it Proof of Theorem \ref{thm:3.2}.}
First of all, we need to show that $|R_5(x)| 
\le 1$ for $x \in T^3$. Solving $\partial_i R_5(x) =0$, $i =1,2,3$, numerically
for critical points shows that $|R_5(x)|$ attains its maximum on the boundary
of $T^3$. Furthermore, 
\begin{align*}
 R_5(x,y,0)& =R_5(x,0,y)=R_5(0,x,y) \\ 
  & = -1+ 2(x + y)-2(x + y)^2 + 
     2(1-4(x + y)+ 4(x^2 + y^2))^2,
\end{align*} 
and the polynomial has no critical point inside $T^2$. Consequently, the 
maximum of $|R_5(x)|$ is attained on the face of $T^3$ defined by $x_1+x_2+x_3
 =1$. In other words, we only need to show that $|U_5(x_1,x_2)|\le 1$ on 
$T^2$. Again this can be proved by solving $\partial_i U_5(x_1,x_2) =0$, 
$i =1,2$, and the maximum is attained on the boundary. This proves that 
$|R_5(x)|\le 1$ on $T^3$ and it also gives the set $\CS_+ = \{x:|R_5(x)| =1\}$
and the set $\CS_- = \{x:|R_5(x)| =-1\}$. Let $S_3$ be the symmetric group
of three elements. For $a = (a_1,a_2,a_3) \in \RR^3$ we define $a\tau := 
(a_{\tau_1},a_{\tau_2},a_{\tau_3})$, $\tau \in S_3$ and 
$(a)_G := \{a\tau: \tau \in S_3\}$. Then 
\begin{align*}
\CS_+ & = 
\{(1/3,1/3,1/3),(0,0,0),(1,0,0)_G,(1/2,1/2,0)_G,(t_1,t_1,1-2t_1)_G\}\\
 \CS_- & = \{((2-\sqrt{2})/4,(2+\sqrt{2})/4,0)_G,(t_2,t_2,1-2t_2)_G\},
\end{align*}  
where $t_1=0.4588164122$ and $t_2= 0.1343303216$. We consider the signature 
$\sigma$ defined by $\sigma(v) = 1$, $v \in \CS_+ \setminus \{(0,0,0)\}$ and  
$\sigma(v) = -1$, $v \in \CS_-$. To show that $\sigma$ is an extremal 
signature, we define $Lf$ by 
\begin{align*} 
   Lf = & c_0 f(1/3,1/3,1/3) + c_1 \sum_{\tau} f((1,0,0)\tau) \\ 
     & + c_2 \sum_{\tau} f((1/2,1/2,0)\tau) + 
     c_3 \sum_{\tau} f((t_1,t_1,1-2t_1)\tau) \\ 
     & - c_4 \sum_{\tau} f(((2-\sqrt{2})/4,(2+\sqrt{2})/4,0)\tau)      
       - c_5 \sum_{\tau} f((t_2,t_2,1-2t_2)\tau),  
\end{align*} 
where the sum is taken over all distinct permutations of the base point and  
the coefficients are given by 
\begin{align*} 
c_0 & =0.0997251873, \quad c_1=0.0097228135, \quad  c_2 =0.0621246411,\\ 
c_3 & = 0.0243979796,\quad c_4=0.0615774830, \quad  c_5= 0.1178707075. 
\end{align*} 
Then $Lf =0$ for all $f \in \Pi_4^3$, which shows that $\sigma$ is an  
extremal signature. This completes the proof of Theorem \ref{thm:3.1}. 
\qed

\medskip

The linear functional $Lf$ given above is evidently a sum of two linear 
functionals with positive coefficients. Unlike the case of $R_3$, however,
the two linear functionals are not cubature formulas of degree $5$ with 
respect to the Lebesgue measure. 

The two cases solved in this section appear to indicate a surprisingly  
complicated picture for the best approximation of monomials in three 
variables, and the picture is remarkably different from that of one and two
variables. We make two remarks in this regard.

\medskip \noindent
{\it Remark 3.1.}
One surprising fact of Theorem \ref{thm:3.2} is that the least deviation is
not given by a reciprocal of an integer. This indicates a major difference 
between the case of three variables and that of one and two variables. In the 
case of one variable, the polynomial of least deviation from zero is the 
classical Chebyshev polynomial, for which the least deviation 
of $x^n$ to $\Pi_{n-1}$ in $C[-1,1]$ is $2^{1-n}$. In the case of two 
variables, we know for example 
$$
 E_n(x^k y^{n-k};B^2)= \inf_{p\in \Pi_{n-1}^2} \|x^k y^{n-k} - p(x)\|_{C(B^2)}
    = 2^{1-n}.
$$
For three variables, however, we do not know if the least deviation of 
$x^\alpha$ to $\Pi_{|\alpha|-1}^d$ could be represented by a simple formula 
that depends only on the total degree of the monomial. The result in this 
section seems to indicate that such a formula does not exist. 
 
\medskip \noindent
{\it Remark 3.2.}
The values of the least deviation in Theorems \ref{thm:3.1} and \ref{thm:3.2}
are surprisingly small. Let us examine the case of the unit ball. We know
\begin{equation}\label{eq:3.1}
 E_n(x_1^k x_2^{n-k};B^3) =  E_n(x_1^k x_3^{n-k};B^3) = 
  E_n(x_2^k x_3^{n-k};B^3) = 2^{1-n}. 
\end{equation}
This follows from the fact that an extremal polynomial $p^*$ for $x_1^m x_2^n$ 
must be even in $x_3$ since $x_1^m x_2^n$ is invariant under the group $\ZZ_2$
applied on the third variable. Let $p^*$ be so chosen; then
\begin{align*}
 \|x_1^{n-m} x_2^m - p^*(x_1,x_2,x_3)\|_{C(B^3)}
 & \ge  \|x_1^{n-m} x_2^m - p^*(x_1,x_2, \sqrt{1-x_1^2-x_2^2})\|_{C(B^2)} \\
& \ge \inf_{p \in \Pi_n^2} \|x_1^{n-m} x_2^m - p(x_1,x_2)\|_{C(B^2)} =2^{1-n}. 
\end{align*}
Furthermore, the equality holds since an extremal polynomial for 
$x_1^{n-m}x_2^m$ on $B^2$ can also serve as an extremal polynomial on $B^3$.
Below is a list of other cases that we know on the unit ball:
\begin{align*}
& E_n(x_1x_2x_3;B^3) = 3^{-3/2}, \quad 
E_n(x_1^2x_2^2x_3^2;B^3) = 2^{-6}\cdot 3^{-2}, \\ 
& E_n(x_1^4x_2^4x_3^4;B^3) =0.5340799374 \cdot 2^{-12}\cdot 3^{-12},
\end{align*}
where we rewrite the value of the third one, which is given in 
Theorem \ref{thm:3.2}, for easier comparision. The value of 
$E_n(x_1^2x_2^2x_3^2;B^3)$ appears to be strikingly small. For other degree 
$6$ monomials given in \eqref{eq:3.1}, the value of the best approximation is 
only $2^{-5}$. Also note the fast decrease shown in these three values. 

\section{Least deviation from zero for $d >3$}
\setcounter{equation}{0} 

We consider the best approximation to $(x_1 \cdots x_d)^2$ on $B^d$ or 
$S^{d-1}$, and the best approximation to $x_1 \cdots x_d$ on $T^d$ or 
$T^{d-1}$ in this section. The extremal polynomial can be taken as symmetric 
polynomials by Proposition \ref{prop:2.3}. It is well known that every 
symmetric polynomial can be written in terms of elementary symmetric 
polynomials (\cite{M}). 

The elementary symmetric polynomials of degree $k$ in variables $x_1, x_2, 
\ldots,x_N$ are defined by 
$$
e_k(x) = \sum_{1\le i_1 < \ldots < i_k \le N}  
  x_{i_1}x_{i_2}\cdots x_{i_k}, \qquad  1 \le k \le N.
$$
In particular, $e_1(x_1,\ldots,x_N) = x_1+ \cdots + x_N$ and 
$e_d(x_1,\ldots,x_N) = x_1 \cdots x_N$. As it is often the case with the
symmetric functions, we assume that $N$ is sufficiently large and do not write 
the dependence of $e_k$ on the number of variables. We will use the notation 
${\bf 1}^k = (1,1,\ldots,1) \in \RR^k$. 

\begin{defn}
Using elementary symmetric functions, define $T_3(x)$ by
$$
T_3(x) = 72 e_3(x) - 4 e_1(x) + 4 e_1^2(x) - 8 e_2(x) +1
$$
and $T_k(x)$ for $k >3$ by the recursive formula
$$
T_k(x) = r_k e_k(x) - T_{k-1}(x), 
$$
where the constant $r_k$ is determined by  $r_k = k^k 
[T_{k-1}(k^{-1}{\bf 1}^k)+1 ]$. 
\end{defn}

Note that $k^{-1}{\bf 1}^k = (k^{-1},\ldots,k^{-1})\in \RR^k$; we use the 
evaluation of $T_{k-1}$, as a function of $\RR^k$, at this point in the 
definition of $r_k$. Clearly $r_k$ is uniquely determined. For $x \in \RR^d$, 
the function $T_d(x)$ will serve as extremal polynomials. In particular,  
the polynomial $T_3(x)$ for $x\in \RR^3$ is the same as $R_3(x)$ in the  
previous section. For $x \in \RR^4$, the explicit formula of $T_4$ is  
given by 
\begin{align*}
T_4(x) = & \, 896x_1x_2x_3x_4 -72(x_1x_2x_3 +x_1x_2x_4+x_1x_3x_4+x_2x_3x_4)\\ 
   &  + 4(x_1+x_2+x_3+x_4) - 4 (x_1+x_2+x_3+x_4)^2 \\ 
  & + 8 (x_1x_2+ x_1x_3+ x_1x_4+ x_2x_3+ x_2x_4+ x_3x_4) - 1.   
\end{align*}
The value of $r_d$ is of particular importance. It can be computed using the 
following formula.

\begin{lem} \label{lem:4.2}
For $d \ge 3$, 
$$
r_d = d \sum_{k=4}^d k^{d-3} \binom{d}{k} \left[(-1)^k (9 k^2 -32 k + 24) + 
 k^2 \right].
$$ 
In particular, $r_3= 72$, $r_4 = 896$, $r_5 = 14400$, and $r_6 =283392$.
\end{lem}

We defer the proof to the end of the section and continue to state our main  
result of this section. 

\begin{thm} \label{thm:4.1} 
For $d \ge 3$, on the $d$-dimensional simplex 
\begin{align*} 
E_{d-1}(x_1\cdots x_d;T^d) & =E_{d-1}(x_1\cdots x_{d-1}(1-x_1-\cdots-x_{d-1}); 
T^{d-1}) \ge r_d^{-1},
\end{align*} 
and the equality holds for $d = 3,4,5$ with $r_d^{-1} T_d(x)$ as a polynomial 
of least deviation from zero. Furthermore, on $B^d$ and $S^{d-1}$,    
$$
E_{2d-1}(x_1^2\cdots x_d^2; B^d) = E_{2d-1}(x_1^2\cdots x_d^2; S^{d-1}) 
\ge r_d^{-2}
$$
and the equality holds with $r_d^{-2}T_d(x_1^2,\ldots,x_d^2)$ as a polynomial
of least deviation from zero.  
\end{thm}

We believe that the equality still holds for $d \ge 6$. In fact, all that is
missing is to prove that $|T_d(x)| \le 1$ for $x \in T^d$. We state it 
as a conjecture.

\begin{con}
For $d \ge 6$, the inequality $\|T_d(x)\|_{C(T^d)} \le 1$ holds. In particular,
the equality in the above two theorems holds for $d \ge 6$.
\end{con}

Let us point out that there does not appear to exist a closed formula for
$r_d$. Below is a list of the first few values of $r_d$ and their prime 
factorization:
\begin{align*}
& r_3 = 72 = 2^3 \cdot 3^2, \quad  r_4 = 896 = 2^7 \cdot 7, \quad 
r_5 = 14400 = 2^6\cdot 3^2\cdot 5^2, \\
& r_6 = 283392 = 2^8\cdot 3^3 \cdot 41, \quad 
r_7 = 6598144 = 2^9\cdot 7^2 \cdot 263, \\
& r_8 = 177373184 = 2^{15} \cdot 5413,\quad  
r_9 = 5406289920 = 2^{12}\cdot 5 \cdot 3^4 \cdot 3259, \\ 
& r_{10} = 184223744000 = 2^{14} \cdot 5^3 \cdot 23 \cdot 3911,\\ 
& r_{11} = 6939874934784 = 2^{14} \cdot 3^3 \cdot 11^2 \cdot 137 \cdot 409. 
\end{align*}
A closed formula will have to catch the pattern of the prime numbers presented 
in these formulae, which seems unlikely. We also note that the values of 
$r_d$ appear to indicate that the best approximation to monomials becomes 
increasingly more complicated as $d$ increases. See also the Remark 3.1.

The proof of Theorem \ref{thm:4.1} is split into several propositions. As
before, we only need to prove the case of polynomials on the simplex.
We start with the point set at which $T_d(x) = \pm 1$. 

We will work with sets of points that are invariant under the symmetric group
$S_d$. For $a \in \RR^d$ we use the notation $(a)_G$ to denote the set of 
points that consist of all distinct permutations of $x$; that is
$$
   (a_1,\ldots,a_d)_G = \{ (a_{\tau_1},\ldots,a_{\tau_d}): 
   \tau = (\tau_1, \ldots, \tau_d) \in S_d\} 
$$   
and we sometimes write $a\tau = (a_{\tau_1},\ldots,a_{\tau_d})$ for 
$\tau \in S_d$.  
 
\begin{prop}
Let $\CS_+(T^d)$ and $\CS_-(T^d)$ be the subsets of $T^d$ on which 
$T_d(x) = 1$ and $T_d(x) = -1$, respectively. For odd $d$,  

\noindent
$\CS_+=\left \{\left(\frac{1}{d}, \ldots, \frac{1}{d}\right)_G, 
 \left( \frac{1}{d-2}, \ldots, \frac{1}{d-2},0,0\right)_G,
\ldots,\left(\frac{1}{3}, \frac{1}{3},\frac{1}{3},0, \ldots, 0\right)_G, 
   (1, 0, \ldots, 0)_G \right \}
$

\noindent
is a subset of $\CS_+(T^d)$ and 

$\CS_-=\left \{ \left(\frac{1}{d-1}, \ldots, \frac{1}{d-1},0\right)_G, 
  \left(\frac{1}{d-3}, \ldots, \frac{1}{d-3},0,0,0\right)_G,
  \cdots, \left(\frac{1}{2}, \frac{1}{2}, 0,\ldots,0\right)_G 
  \right \}
$

\noindent
is a subset of $\CS_-(T^d)$. For even $d$,  

$\CS_+=\left \{\left(\frac{1}{d}, \ldots, \frac{1}{d}\right)_G, 
 \left( \frac{1}{d-2}, \ldots, \frac{1}{d-2},0,0\right)_G,
  \cdots, \left(\frac{1}{2}, \frac{1}{2}, 0,\ldots,0\right)_G 
  \right \}
$

\noindent
is a subset of $\CS_+(T^d)$ and 

$\CS_-=\left \{ \left(\frac{1}{d-1}, \ldots, \frac{1}{d-1},0\right)_G, 
  \left(\frac{1}{d-3}, \ldots, \frac{1}{d-3},0,0,0\right)_G,
\ldots,   (1, 0, \ldots, 0)_G \right \}
$

\noindent
is a subset of $\CS_-(T^d)$. Furthermore, all these points are on the 
face of $T^d$ defined by the equation $x_1 +\ldots +x_d =1$. 
\end{prop}

\begin{proof}
In the definition of $T_d$, the value of $r_d$ is chosen so that 
$T_d(\frac{1}{d+1}, \ldots, \frac{1}{d+1}) =1$. All other points in the
given set contain at least one zero component. This allows us to use
induction. For $T_3(x)$, it is easy to verify that $T_3(x) =1$ if 
$x = (1,0,0)_G$ and $x =(1/3,1/3,1/3)$, and $T_3(x) =-1$ if $x=(1/2,1/2,0)_G$.
The induction is based on the formula
$$
  T_d(x_1,\ldots, x_{d-1},0) = - T_{d-1}(x_1,\ldots,x_{d-1})
$$
and similar formulae obtained by a permutation of $(x_1,\ldots,x_{d-1})$, 
which follow from the definition of $T_d$ and the fact that 
$e_d(x_1,\ldots,x_{d-1},0)=0$.
\end{proof}

\begin{prop}
The signature $\sigma$, defined by $\sigma(v) =1$ if $v\in \CS_+$ and 
$\sigma(v) =-1$ if $v\in \CS_-$, is an extremal signature of $T_d$. More
precisely, define 
\begin{align*}
 Lg &=  d^{d-1} g\left(\frac{1}{d}, \ldots, \frac{1}{d}\right) - (d-1)^{d-1} 
\sum_\tau g\left(\left(\frac{1}{d-1},\ldots,\frac{1}{d-1},0\right)
  \tau\right)  \\
& + \ldots + (-1)^{d-2} 2^{d-1} 
\sum_{\tau} g\left(\left( \frac{1}{2},\frac{1}{2},
    0,\cdots,0\right)\tau \right)+ (-1)^{d-1}\sum_{\tau} 
   g((1, 0, \ldots, 0)\tau);
\end{align*}
then $L p =0$ for all $p \in \Pi_{d-1}^d$.
\end{prop}

\begin{proof}
Since the points in $\CS_+$ and $\CS_-$ are symmetric with respect to $S_d$, 
a moments reflection shows that we only need to verify $L g=0$ for symmetric 
polynomials. One basis of symmetric polynomials in $\Pi_{d-1}^d$ consists of
$m_k(x) = x_1^k +\cdots + x_d^k$ for $k =0, 1,\ldots,d-1$. We show $L m_k =0$.
Let $a_j = (1/j,\ldots,1/j,0,\ldots,0) \in T^d$, which contains exactly $j$ 
nonzero entries. Then the sum $\sum_{\tau\in S_d} g(a_j\tau)$ contains 
$\binom{d}{j}$ terms and $m_k(a_j) = j (1/j)^k$. Consequently, for $k \ge 1$, 
$$
 L m_k = d^{d-k} - \binom{d}{1} (d-1)^{d-k} + \binom{d}{2} (d-2)^{d-k} + 
   \ldots + (-1)^{d-1} \binom{d}{d-1}. 
$$
Furthermore, it is easy to see that $L 1$ gives the same formula as $L m_1$ 
(recall that points in $\CS_+$ and $\CS_-$ satisfy $x_1+\ldots + x_d =1$).
Hence, we need to show, changing $n-k$ to $k$,   
\begin{equation} \label{eq:combi}
  \sum_{j=0}^d (-1)^j \binom{d}{j} j^k =0, \qquad k = 1,2,\ldots,d-1.
\end{equation}
This is well known and can be proved by induction on $d$. 
\end{proof}

By Theorem 3.2, the above proposition has proved that 
$$
 \inf_{p \in \Pi_{d-1}^d} \|x_1\ldots x_d - p(x)\|_{C(T^d)} \ge r_d^{-1}.
$$
In order to show that the equality holds, we only need to prove that 
$|T_d(x)| \le 1$ for $x \in T^d$. However, we are able to establish this 
inequality only for $d =3,4,5$. 

\begin{prop}
For $d =3, 4, 5$, 
$$
  |T_d(x)| \le  \|T_d\|_{C(T^{d-1})} =1, \qquad x \in T^d.  
$$
\end{prop}

\begin{proof}
There does not seem to be an easy way of proving this. We use the standard
method of finding critical points upon solving $(\partial_i T_d(x)/
\partial x_i) = 0$, $1\le i \le d$. The case $d=3$ is in the previous section
and the equations can be solved algebraically. The cases $d =4$ and $d=5$ are 
solved numerically. The details are omitted. Once the critical points are
found, we can then verify that the inequality $|T_d(x)| < 1$ holds on these
points, which shows that the maximum of $|T_d(x)|$ is attained on the 
boundary of $T^d$. Since $T_d(x_1,\ldots,x_{d-1},0) = - T_{d-1}(x_1,\ldots,
x_{d-1})$, by induction, we only need to prove that $|T_d(x_1,\ldots,x_{d-1},
1-x_1-\ldots -x_{d-1})| \le 1$ for $x \in T^{d-1}$. Again, this is done by 
computing the critical points and evaluating. 
\end{proof}

Putting the above propositions together, we have completed the proof of 
Theorem \ref{thm:4.1}. 

We still need to prove Lemma \ref{lem:4.2}. First we note that the definition 
of $T_d$ implies
$$
  T_d(x) = \sum_{k=4}^d (-1)^{d-k} r_k e_k(x) + (-1)^{d-3}T_3(x).
$$

\medskip\noindent
{\it Proof of Lemma \ref{lem:4.2}.}
Setting $x = a_{d+1}= (d+1)^{-1} {\bf 1}^{d+1}$ and using the fact that 
$e_k(a_{d+1}) = \binom{d+1}{k} (d+1)^{-k}$ leads to the relation
$$
(d+1)^{-(d+1)} r_{d+1} = \sum_{k=4}^{d}(-1)^{d-k} \binom{d+1}{k} (d+1)^{-k}r_k 
   + (-1)^{d-3} T_3(a_{d+1})+1. 
$$
Replacing $d+1$ by $d$ we can write the above equation as 
$$
\sum_{k=4}^d (-1)^{d-k} \binom{d}{k} d^{-k}r_k = (-1)^d T_3(d^{-1} {\bf 1}^d) 
  +1:= A_d.
$$
We want to reverse this relation so that $r_d = \sum_{j=4}^d b_j A_j$. A 
standard argument shows that, with $b_j$ given below, this will follow from
the combinatoric relation
$$
J_{k,d}:= \sum_{j=k}^d b_j (-1)^{j-k}\binom{j}{k}j^{-k}=\delta_{k,d}, \qquad 
b_j = d \binom{d}{j} j^{d-1},
$$
where $\delta_{k,d} =1$ if $k=d$ and $0$ otherwise. For $k =d$, $J_{d,d} =1$
holds trivially. For $k < d$, a change of summation index gives 
$$
J_{k,d} = \binom{d}{k} \sum_{j=k}^{d} (-1)^{j-k}\binom{d-k}{j-k} j^{d-1-k}
= \binom{d}{k} \sum_{j=0}^{d-k} (-1)^{j}\binom{d-k}{j}(k+j)^{d-1-k}.
$$
By \eqref{eq:combi}, the last formula shows that $J_{k,d}=0$. Finally, the 
definition of $T_3$ shows that 
$$
T_3(d^{-1} {\bf 1}^d) = 72 \binom{d}{3} d^{-3} - 8\binom{d}{2}d^{-2} + 1
  = d^{-2}(9 d^2 - 32 d +24),
$$
which gives the explicit value of $A_d$. Putting these relations together
completes the proof.
\qed

\medskip

We conjecture that $|T_d(x)|\le 1$ for all $d$ and the points in $\CS_+$ and
$\CS_-$ are the only ones on the face of $T^d$ defined by $x_1+\ldots+x_d =1$
on which $|T_d(x)|=1$. The following fact is helpful for the case $d =5$ and 
could be useful for $d >5$. Let 
$$
D_d(x) = \det \left[\begin{matrix} 1&1& \ldots &1\\x_1&x_2& \ldots &x_d\\
    \vdots& \vdots& \ddots &\vdots \\ x_1^{d-3}&x_2^{d-3}& \ldots &x_d^{d-3}\\
    \partial_1 T_d(x) &  \partial_1 T_2(x) & \ldots & \partial_d T_d(x)   
    \end{matrix}\right],
$$
where $\partial_i = \partial/\partial x_i$. Then $D_d(x)$ can be factored 
completely. We have, for example, 
$$
 D_5(x) = -64 \prod_{1\le i < j \le 4}(x_i - x_j) (-14 + 225 x_5). 
$$
At the critical points of $T_5$, $D_5(x) = 0$ and $D_5(x\tau) =0$ for 
$\tau \in S_5$. One could use it to confirm the conjecture for $d > 5$. We 
did not try hard to push for larger $d$, since the method does not seem to 
lead to a proof for all $d$. 

Our conjecture implies that $T_d(x)$ attains its maximum on the boundary of 
$T^d$. Part of this can be proved as follows: Let $\Delta$ denote the 
Laplacian operator $\Delta = \partial_1^2 + \ldots + \partial_d^2$. Then it
is easy to verify that $\Delta T_d(x) = (-1)^{d-1} 8$. In particular, this
shows that $(-1)^{d-1} T_d(x)$ is a subharmonic function. Hence, by the maximum
principle (\cite{J}) for the subharmonic functions, we can conclude that 
$(-1)^{d-1} T_d(x) \le \max_{x \in \partial T^d} (-1)^{d-1} T_d(x)$.

\bigskip
{\it Acknowlegement.} The author thanks the referees for their careful review.

\enddocument